\documentclass[12pt,reqno]{amsart}

\newcommand{\R}{{\mathbb R}}

\newcommand{\e}{\epsilon}

\newcommand{\D}{\Delta}
\newcommand{\p}{\partial}

\numberwithin{equation}{section}
\newtheorem{theorem}{Theorem}[section]
\newtheorem{defn}[theorem]{Definition}
\newtheorem{lemma}[theorem]{Lemma}
\newtheorem{remark}[theorem]{Remark}

\newtheorem{coro}[theorem]{Corollary}

\begin{document}

\title[SPDEs under fast   dynamical boundary conditions ]
{Reductions and deviations for stochastic partial differential
equations under fast dynamical boundary conditions }

\author[W. Wang \&  J. Duan  ]
{Wei Wang and Jinqiao Duan }

\address[W.~Wang ]
{Department of Mathematics\\ Nanjing University\\
Nanjing, 210093, China } \email[W.~Wang]{wangweinju@yahoo.com.cn}

\address[J.~Duan]
{Department of Applied Mathematics\\
Illinois Institute of Technology\\
Chicago, IL 60616, USA} \email[J.~Duan]{duan@iit.edu}

\date{August 7, 2008}

\thanks{This work was partly supported by the NSFC grant
   10701072 and the NSF grant 0620539.}



\begin{abstract}
In order to understand the impact of
  random influences at physical boundary on the evolution of
   multiscale systems, a  stochastic partial
differential equation model under a fast random dynamical boundary
condition is investigated. The noises in the model and in the
boundary condition  are both \emph{additive}. An effective equation
is derived and justified by reducing the random \emph{dynamical}
boundary condition  to a simpler one. The effective system is still
a stochastic partial differential equation. Furthermore, the
quantitative comparison between the solution of the original
stochastic system and the effective solution is provided by
establishing normal deviations and large deviations principles.
Namely, the normal deviations are asymptotically characterized,
while the rate and speed of the large deviations are estimated.

\end{abstract}

\maketitle

\section{Introduction}\label{s1}

The    random  fluctuations may have delicate impact in modeling,
analyzing, simulating and predicting complex phenomena. The need to
quantify uncertainties has been widely recognized in geophysical and
climate dynamics, materials science, chemistry, biology and other
areas \cite{Imkeller, E00,TM05}. Stochastic partial differential
equations (SPDEs or stochastic PDEs) are appropriate mathematical
models for various  multiscale systems under random influences
\cite{WaymireDuan}.

A stochastic partial differential equation usually contains noise
in the equation itself, i.e., the noise is acting on the system
inside the physical medium \cite{Roz, PZ92, Chow, Walsh}.
 However, noise may
affect  a complex system not only inside the physical medium but
also on the physical boundary. Randomness in such boundary
conditions are often due to various fast time scale environmental
fluctuations.

 The usual  boundary conditions, such as the   Dirichlet or Neumann
boundary conditions, do not contain time derivatives of the system
state. On the contrary, \emph{dynamical} boundary conditions
contain time derivatives of the state.

 The   boundary conditions may further contain
random effects, as in some applications. For example, the
environment surrounding a pipe fluid is usually subject to
uncertain fluctuations, such as random vibration around a natural
gas pipe or a waste water pipe. In a fluid laboratory, a wind
tunnel or fluid pipe  may be sitting on a flat foundation, which
is also subject to random vibration. This noise affects the pipe
fluid flow via boundary conditions such as a Dirichlet or Neumann
boundary condition  on a part of the pipe surface boundary, which
is   a static  boundary condition perturbed by random
fluctuations. The salinity flux on fluid inlet boundary of a
gravity current (e.g., at the Strait of Gibraltar) has a
fluctuating component and this leads to a random Neumann boundary
condition \cite{DijkBook, Bong}. Stochastic dynamical systems
under such random boundary conditions have been studied recently
in, for example, \cite{DaPrato3, FW92, Mas95, Sowers, DuanGaoSchm,
Bong}.


 In some other   applications, the evolution of systems
 may also be subject to   dynamical boundary conditions (containing time
derivatives of the system state), under random perturbations. Such
\emph{random  dynamical  boundary conditions} arise in the
modeling of, for example, the air-sea interactions on the ocean
surface \cite{PeiOor92}, heat transfer in a solid in contact with
a fluid \cite{Langer}, chemical reactor theory \cite{Lap}, as well
as colloid and interface chemistry \cite{Vold}. The random
fluctuations on the boundary are usually much faster than the
intrinsic time scale of these systems \cite{DuanGaoSchm}.
 In these cases, the mathematical models are  stochastic PDEs with fast-varying random
dynamical boundary conditions. Such stochastic dynamical systems
have been more recently investigated in,
for example, \cite{DuanGaoSchm, ChSch04, ChSch05, YD05, Bona}.\\

Motivated by better analytical understanding of the
above-mentioned multiscale systems under fast scale random
dynamical interactions on the physical boundary, as well as noisy
forcing inside the physical medium,   we consider a stochastic
parabolic partial differential equation on a bounded domain $D$
under fast varying random dynamical boundary condition on a part
of the boundary $\p D$. The fast time scale in the random
dynamical boundary condition is controlled by a small parameter
$\e>0$ and perturbed by a   noise (white in time but correlated in
space). Specifically, we study the following  stochastic parabolic
partial differential equation with a random dynamical boundary
condition on a part of boundary $\Gamma_1$ and a Dirichlet
boundary condition on the rest of boundary:
\begin{eqnarray*}
du_\e&=&[\D u_\e+f(u_\e)]\,dt+\sigma_1 \,dW_1(t)\,,\;\; in \;\;D\\
\e du_\e&=&[-\p_\nu u_\e-u_\e]\,dt +\sqrt{\e}
\;\sigma_2\,dW_2(t)\,,\;\;
on\;\; \Gamma_1 \\
u_\e&=& 0\,, \;\; on \;\; \Gamma_2
\end{eqnarray*}
where $\Gamma_1$ and $\Gamma_2$ form the whole  of the boundary $\p
D$ of domain $D$, $f(u)$ is some nonlinear term, and $\sigma_1,
\sigma_2$ and $\e$ are constants. More details of
this model will be presented in the next section.\\

 First,  we derive an effective    model
 for the above system as $\e\rightarrow 0$\,. The effective model
 is still a stochastic partial differential equation,
 but with a simpler  boundary condition
 (see  Theorem \ref{main1}). Note that simple boundary
 conditions not only facilitate    theoretical
 analysis but also are desirable for numerical simulations.
  To this end, we formulate the above stochastic system as
  an abstract
 stochastic evolution equation with non-Lipschitz nonlinear term having
 polynomial growth and fulfilling a suitable dissipativity condition.
 Since the nonlinear term
 is non-Lipschitz, we introduce a cut-off function and a stopping time to
 obtain a unique mild solution, which is also the unique weak
 solution for the system (see \S 3). Here we follow a semigroup approach  which
 is also used in \cite{Cerrai} to obtain a mild solution for a class of
 reaction-diffusion equations with  multiplicative noise and  non-Lipschitz
 reaction term but with deterministic static homogenous boundary condition.
Next we present some useful a priori estimates for the weak
solutions which yield the tightness of the distributions of the
solutions (see \S 4). Then by a discussion in the variational form
of the system we can pass the limit $\e\rightarrow 0$ to get the
effective equation which is a stochastic partial differential
equation but with a simpler boundary condition. And further we show
that $u_\e$ converges to the effective solution $u$ of the limiting
equation, in probability, in an appropriate function space.

Then we determine the normalized deviation between the solution
$u_\e$ of the original system and the solution $u$ of the effective
system we obtained. it is proved that the normalized deviation
\begin{equation*}
v_\e=\frac{1}{\sqrt{\e}}[u_\e-u]
\end{equation*}
converges, as $\e\rightarrow 0$, to a  process which solves
 a linear   partial differential equation   with random
 coefficients, under a
 \emph{random} boundary condition. Namely, this random
  boundary condition is a deterministic
 static  boundary condition perturbed by a white noise; see
 Theorem \ref{main2}.

 Finally, we investigate the  deviations $u_\e-u$ of
 order $\e^\kappa$, with $0<\kappa<1/2$.
 In fact a large deviation result (see Theorem \ref{main3}) is proved
for

\begin{equation*}
\frac{1}{\e^\kappa}[u_\e-u].
\end{equation*}
In this weak convergence approach we prove the Laplace principle
which is equivalent to the large deviations principle in a complete
separable metric space (i.e., Polish space).  For background see
\cite{BD00, Dupuis,Feng}.

In \cite{FW92} the authors have studied a system of
reaction-diffusion equations with Lipschitz nonlinear term in open
interval $(-1, 1)$ with a static boundary condition perturbed by a
  stationary random process which varies fast in time. The
solution is represented through a Green function, then the limit
$\e\rightarrow 0$ is passed in the space $C\big([0,T]\times (-1,
1)\big)$ and one gets a deterministic partial differential
equation. The normal deviations of the solutions are then obtained
in a weighted space
with some  assumptions on mixing properties of the random stationary process.\\

This paper is organized as follows. The problem formulation and some
preliminary results are presented in \S 2\,. Section 3 is devoted to
the derivation of some useful a priori estimates. The tightness of
the laws or distributions for the solutions is proved in \S 4 and
then the   effective model is derived in \S 5\,. The normalized
deviations are studied in section \S 6 and the last section, \S 7\,,
is devoted to a large deviations result.

\vskip 0.8cm


\section{Problem formulation}

Let $D$  be a bounded  smooth domain in $\R^N$ ($1\leq N \leq 3$),
with boundary $\Gamma_0$\,. Assume that
$\Gamma_0=\overline{\Gamma}_1\cup\overline{\Gamma}_2$\,, where
$\Gamma_1$ and $\Gamma_2$ are open subsets  of $\Gamma_0$ and
$\Gamma_1\cap\Gamma_2=\emptyset$\,. We consider the following
stochastic partial differential equation  with a random dynamical
boundary condition on $\Gamma_1$ and a Dirichlet boundary condition
on $\Gamma_2$\,:
\begin{eqnarray}\label{fs1}
du_\e&=&\big [\D u_\e +f(u_\e) \big]\,dt+\sigma_1\, dW_1(t)\;\;in\;\; D\,,\\
\e du_\e&=&\big[ -\p_\nu u_\e-u_\e\big]\,dt+\sqrt{\e} \,\sigma_2 \,dW_2(t)\;\; on \;\; \Gamma_1\,, \label{fs3}\\
u_\e&=&0\;\; on\;\; \Gamma_2\,,\\
u_\e(0)&=&u^0\;\; in\;\; D\,, \\
\gamma_1u_\e(0)&=&\gamma_1u^0\;\; on \;\;\Gamma_1 \label{fs5}
\end{eqnarray}
where $\e$ is real number with $0<\e\leq 1$\,; $\nu=(n_1, n_2,
\cdots, n_N)$ is the outer unit normal vector on $\Gamma_1$\,, and
$$
\p_\nu=\sum^N_{i=1}n_i\p_{x_i}\,.
$$
$W_1(t,x)$ and $W_2(t,x)$ are mutually independent $L^2(D)$-valued
and $L^2(\Gamma_1)$-valued Wiener processes, respectively, on a
complete probability space $(\mathbf{\Omega}, \mathcal{F},
\mathbf{P})$ with a canonical filtration $(\mathcal{F}_t)_{t\geq
0}$\,. Moreover, $\gamma_1$ is a trace operator on $\Gamma_1$ (see
next paragraph), and $\sigma_1$ and $\sigma_2$ are both constants.
In fact, $\sigma_1:=\sigma_1\, I_{L^2(D)} $ and $\sigma_2:=\sigma_2
\; I_{L^2(\Gamma_1)}$\,, where $I_{L^2(D)}$ and $I_{L^2(\Gamma_1)}$
are the identity operators on $L^2(D)$ and $L^2(\Gamma_1)$
respectively. Taking $\e=0$ we formally have
\begin{eqnarray}\label{avsys1}
du&=&\big [\D u +f(u) \big]\,dt+
\sigma_1\,dW_1(t)\;\;in\;\; D,\\
\p_\nu u+u&=&0\;\; on \;\; \Gamma_1\,, \label{avsys3}\\
u&=&0\;\; on\;\; \Gamma_2\,,\\
u(0)&=&u^0\;\; in\;\; D\,. \label{avsys4}
\end{eqnarray}

In the following we rewrite the equation in an abstract setting.
Denote by $\mathbf{H}^s(D)$ the Sobolev spaces $W_2^s(D)$\,, $s\geq
0$\,, with the usual norms; see for instance \cite{Lion1} for the
definition. Note that $\mathbf{H}^0(D)=L^2(D)$\,. We also define
$\mathbf{H}^s_{\Gamma_2}(D)$ and $\mathbf{H}_0^s(D)$\,, $s\geq 0$\,,
the spaces consisting of the function of $\mathbf{H}^s(D)$ which
vanish on the boundary $\Gamma_2$ and $\p D$\,, respectively. In
addition we denote by $\mathbf{B}^s(\Gamma_1)$ the Sobolev space
$\mathbf{H}^{s-\frac{1}{2}}(\Gamma_1)$\,, $s>\frac{1}{2}$ on the
boundary $\Gamma_1$\,. And let $\gamma_1$ be the trace operator with
respect to $\Gamma_1$ which is continuous linear operator from
$\mathbf{H}^s(D)$ to $\mathbf{B}^s(\Gamma_1)$ for $s>\frac{1}{2}$\,.
For more information on trace operators we refer to \cite{Tri78}. We
also define spaces $\mathbf{H}^{-s}(D)$\,, $s\geq 0$ and
$\mathbf{B}^{1-s}(\Gamma_1)$, $s>\frac{1}{2}$ as the dual spaces of
$\mathbf{H}_0^s(D)$ and $\mathbf{B}^{s}(\Gamma_1)$ respectively.

We denote by $\langle\cdot\,, \cdot\rangle_D$\,, $\langle\cdot\,,
\cdot\rangle_{\Gamma_1}$ the usual inner products in $L^2(D)$ and
$L^2(\Gamma_1)$\,, respectively.

For our system we introduce the following functional spaces
$$
X^1_\e=\Big\{ (u, v )\in
\mathbf{H}_{\Gamma_2}^1(D)\times\mathbf{B}^1(\Gamma_1):
v=\sqrt{\e}\gamma_1 u \Big\}
$$
and $X^0=L_{\Gamma_2}^2(D)\times L^2(\Gamma_1)$\,. Here
$L_{\Gamma_2}^2(D)$ is the space consisting of functions in $L^2(D)$
which vanish on $\Gamma_2$\,. Define the norm and inner product on
$X^0$ and $X^1_\e$ respectively as
$$|(u, v)|^2_{X^0}=|u|^2_{L^2(D)}+|v|^2_{L^2(\Gamma_1)}\,,\;\; (u, v)\in X^0\,,$$
$$\big\langle (u_1, v_1) , (u_2, v_2)\big\rangle_{X^0}=\langle u_1, u_2\rangle_D+\langle v_1, v_2\rangle_{\Gamma_1}$$
and
$$|(u, v)|^2_{X_\e^1}=|u|^2_{\mathbf{H}^1_{\Gamma_2}(D)}+|v|^2_{\mathbf{B}^1(\Gamma_1)}\,,\;\; (u, v)\in X_\e^1\,,$$
$$\big\langle (u_1, v_1)\,, (u_2, v_2)\big\rangle_{X^1_\e}=\langle u_1\,, u_2\rangle_{\mathbf{H}_{\Gamma_2}^1(D)}+
\langle v_1\,, v_2\rangle_{\mathbf{B}^1(\Gamma_1)}$$ for $(u_i\,,
v_i)\in X^0$\,, $i=1, 2$\,. Here $|u|_{\mathbf{H}^1_{\Gamma_2}(D)}$
is taken as the equivalent norm as
$|u|_{\mathbf{H}^1_{\Gamma_2}(D)}=|\nabla u|_{L^2}$\,.

\begin{remark}
Notice that the space $X_\e^1$ depends on $\e$ but in the next
section we just give some estimates for a fixed $\e>0$\,. And for
passing the limit $\e\rightarrow 0$ we will consider in a fixed
space instead of $X_\e^1$\,.
\end{remark}

Now we define a boundary operator $B$ on $\mathbf{H}^1(D)$ as
\begin{equation*}
Bu=\gamma_1\p_\nu u+\gamma_1u\,,\;\; u\in \mathbf{H}^1(D)
\end{equation*}
and a second order differential operator $A=-\D$ with homogenous
Neumanna boundary condition. Then we introduce the operator
$\mathcal{A}$ on $D(\mathcal{A})=\{(u\,, v)\in X^1_\e: (-\D u, B
u)\in X^0 \}$ as
\begin{equation}\label{A}
\mathcal{A}z=\Big(Au\,, \frac{1}{\sqrt{\e}}Bu\Big )\,, \;\; z=(u\,,
v)\in X^1_\e\,.
\end{equation}
Associate with the operator $\mathcal{A}$ we introduce the following
bilinear form on $X^1_\e$
\begin{equation}\label{bi}
a(z, \bar{z})=\big\langle\mathcal{A}z, \bar{z}
\big\rangle_{X^0}=\int_D\nabla u\nabla\bar{u}\,dx
+\int_{\Gamma_1}(\gamma_1u)(\gamma_1\bar{u})\,d\Gamma_1
\end{equation}
with $z=(u\,, v)$\,, $\bar{z}=(\bar{u},\bar{v})\in X_\e^1$\,.
Noticing that
\begin{equation}\label{e:1}
|\gamma_1u|^2_{\mathbf{B}^1(\Gamma_1)}\leq
C(\Gamma_1)|u|^2_{\mathbf{H}^1(D)}\,,
\end{equation}
 there is some constants $M>0$\,, $\tilde{\alpha}>0$ and
$\tilde{\beta}\in\R$ such that
\begin{equation}\label{e:a}
a(z,\bar{z})\leq M|u|_{\mathbf{H}^1(D)}|\bar{u}|_{\mathbf{H}^1(D)}
\end{equation}
and
\begin{equation*}
a(z,z)\geq
\tilde{\alpha}|u|^2_{\mathbf{H}^1(D)}-\tilde{\beta}|u|^2_{L^2(D)}\,.
\end{equation*}
Also by (\ref{e:1}) then the following coercive property of $a$
holds
\begin{equation}\label{coercive}
a(z,z)\geq \alpha|z|^2_{X^1_\e}-\beta|z|^2_{X^0}\,,\;\; z\in X^1_\e
\end{equation}
for some constants $\alpha>0$ and $\beta\in\R$. Then the linear
operator $-\mathcal{A}$ generates a $C_0$--semigroup,
$\mathcal{S}(t)$, which is compact and analytic on $X^0$\,; see
\cite{AmEsc96}.

\begin{remark}
Sometimes we also use the notation $a(u, \bar{u})$ instead of $a(z,
\bar{z})$ for any $z=(u, v)$, $\bar{z}=(\bar{u}, \bar{v})\in
X_\e^1$\,.\\
\end{remark}

Now for the nonlinear term we make the following assumptions

\begin{description}
  \item[(F)] $f:\R\rightarrow \R$ is $C^1$-continuous and there are positive constants $a_1$\,, $b_1$ such that
             $$
            f(u)u\leq -a_1u^4+b_1\,,\;\; {\rm for} \;u\in\R
              $$
             $$
            |f(u)|\leq a_1|u|^3+b_1\,,\;\; {\rm for} \; u\in \R
              $$
              $$
              f'(u)\leq -a_1u^2+b_1\,,\;\; {\rm for} \; u\in \R\,.
              $$
\end{description}

\medskip

Define $F(u)=\int_0^uf(v)dv$\,. Then by the assumptions
($\mathbf{F}$) there is some positive constant $\tilde{b}$ such that
\begin{equation}\label{FG}
F(u)\leq \tilde{b}\;\; {\rm and} \;f'(u) \leq \tilde{b}\,,\;\;
\forall u\in \R\,.
\end{equation}
For the stochastic term we assume the following conditions.
\begin{description}
  \item[($\Sigma$)] Stochastic process $W(t)=(W_1(t)\,, W_2(t))^t$\,, is a $Q$-Wiener process on $X^0$\,, defined on a
  filtered probability space $(\mathbf{\Omega}\,, \mathcal{F}\,, \mathcal{F}_t\,,\mathbf{P})$ with covariance
  operator $Q=(Q_1\,, Q_2)$ is trace class. Furthermore we assume Tr$(A^{\frac{1}{2}}Q_1)<\infty$\,.

\end{description}

\begin{remark}
An example of such functions $f$  is given by the following cubic
polynomial
$$
f(u)=-au^3+bu^2+c
$$
with $a>0$\,, $b$\,, $c\in\R$\,.
\end{remark}
\begin{remark}
In $(\Sigma)$, the technical condition on $Q_1$ is for the proof of
regular properties of solution. As one example for such $Q_1$\,, one
can define $W_1=\sqrt{q(x)}w(t)$ with positive function $q\in
\mathbf{H}^1_0(D)$ and $w$  a standard scalar Wiener process. Then,
by the property of trace operator, the covariance operator of
$W_1$\,, $Q_1=q$ satisfies
$Tr(A^{1/2}q)=|q|^2_{\mathbf{H}^1_0(D)}<\infty$\,.
\end{remark}

\smallskip

With the above notations system (\ref{fs1})--(\ref{fs5}) can be
written as the following abstract stochastic evolutionary equation
\begin{eqnarray}\label{SEE}
dz_\e(t)=\big[-\mathcal{A}z_\e+H_\e(z_\e)\big]\,dt+\Sigma\,
dW(t)\,,\;\; z(0)=z^0
\end{eqnarray}
where $z_\e=(u_\e\,, \sqrt{\e}\gamma_1u_\e)^t$,
$H_\e(z)(x)=\big(f(u(x)), 0 )\big)^t$,
$\Sigma=\big(\sigma_1,\sigma_2\big)^t$ and $W(t)=(W_1(t)\,,
W_2(t))^t$\,. The equation (\ref{SEE}) can be further rewritten in
the following mild sense
\begin{equation}\label{mild}
z_\e(t)=\mathcal{S}(t)z^0+\int_0^t\mathcal{S}(t-s)H_\e(z_\e(s))\,ds
+\int_0^t\mathcal{S}(t-s)\Sigma\,dW(s)\,.
\end{equation}
An adapted process $z_\e$ is called a mild solution of (\ref{SEE})
if (\ref{mild}) hold.  For $u_\e\in L^2(0, T;
\mathbf{H}_{\Gamma_2}^1(D)) \cap C(0, T; L^2(D))$\,, we call $u_\e$
a weak solution of (\ref{fs1})--(\ref{fs5}) if for any $t\in[0, T)$
\begin{eqnarray}\label{weak}
&&\langle u_\e(t), \psi(t)\rangle_D+\e\langle
\gamma_1u_\e(t),\psi(t) \rangle_{\Gamma_1}\\
&=& \langle u_\e(0),\psi(0)\rangle_D+\e\langle
\gamma_1u_\e(0),\psi(0) \rangle_{\Gamma_1}+\int_0^t\big\langle
u_\e(s), \frac{\p\psi(s)}{\p
t}\big\rangle_D\,ds\nonumber\\&&{}+\e\int_0^t\big\langle
(\gamma_1u_\e)(s),\frac{\p\psi(s)}{\p t} \big\rangle_{\Gamma_1}\,ds
+\int_0^t a(u_\e, \psi)\,ds +\int_0^t\langle f(u_\e),
\psi\rangle_D\,ds\nonumber\\
 &&{}+ \int_0^t\langle\psi,
\sigma_1\,d{W}_1(s)
 \rangle_D+ \sqrt{\e}\int_0^t\langle\psi, \sigma_2\,d{W}_2(s)
\rangle_{\Gamma_1}\,, \nonumber
\end{eqnarray}
for any $\psi\in C^1(0, T; C^\infty(D))$.
For more detail about solution of SPDEs we refer to \cite{PZ92}.

We end this section by recalling the following two lemmas from
\cite{Lions}, which will be used in our later analysis.
\begin{lemma}\label{Lions}
Let $\mathcal{Q}$ be a bounded region in $\R\times \R^n$. For any
given functions $h_\e$ and $h$ in $L^p(\mathcal{Q})$ $(1<p<\infty)$,
if
$$
|h_\e|_{L^p(\mathcal{Q})}\leq C\,,\;\; h_\e\rightarrow h\;\; {\rm
in} \;\;\mathcal{Q}\;\; {\rm almost\; everywhere}
$$
for some positive constant $C$\,, then $h_\e \rightharpoonup h$
weakly in $L^p(\mathcal{Q})$\,.
\end{lemma}

Let $\mathcal{X}\subset \mathcal{Y}\subset \mathcal{Z}$ be three
reflexive Banach spaces and $\mathcal{X}\subset \mathcal{Y}$ with
compact and dense embedding. Define Banach space
$$
\mathcal{G}=\{h: h\in L^2(0 ,T; \mathcal{X}), \frac{dh}{dt}\in
L^2(0, T; \mathcal{Z}) \}
$$
with norm
$$
|h|^2_\mathcal{G}=\int_0^T|h(s)|^2_\mathcal{X}ds+\int_0^T\Big|\frac{dh}{ds}(s)\Big|^2_\mathcal{Z}ds,\;\;
h\in \mathcal{G}.
$$

\begin{lemma}\label{comp}
If $G$ is bounded in $\mathcal{G}$, then $G$ is precompact in
$L^2(0, T; \mathcal{Y})$.
\end{lemma}

\vskip 0.8cm


\section{Some a priori estimates}

In this section we prove the stochastic evolutionary equation
(\ref{SEE}) is well-posed and further derive a few useful a priori
estimates on the solutions. Since the nonlinear term is
non-Lipschitz, we apply the cut-off technique with a random
stopping time. The same idea was used in \cite{Cerrai, Cerrai2}
for stochastic reaction-diffusion equations with local Lipschitz
nonlinear terms and multiplicative noise. See also \cite{LS06} for
systems on unbounded domain.

\begin{theorem}\label{wellpose} (\textbf{Wellposedness})\\
Assume that $(\bf{F})$ and $(\mathbf{\Sigma})$ hold.  For any $T>0$,
let $z^0=(u^0, v^0)\in X^0$ be a $\big(\mathcal{F}_0,
\mathcal{B}(X^0)\big)$-measurable random variable. Then system
(\ref{SEE}) has a unique mild solution $z_\e\in
L^2\big(\mathbf{\Omega}, C(0,T; X^0)\cap L^2(0, T; X_\e^1)\big)$,
which is also a weak solution in the following sense
\begin{eqnarray}
&&\hspace{-0.3cm}\langle z_\e(t),\phi(t)\rangle_{X^0}-\langle
z_\e(0),\phi(0)\rangle_{X^0} \label{weak1}\\
\hspace{-0.3cm}&=&\hspace{-0.3cm}\int_0^t\big\langle z_\e(s),
\frac{\p \phi(s)}{\p t}\big \rangle_{X^0}\,ds
-\int_0^t\langle\mathcal{A} z_\e(s),
\phi(s)\rangle_{X^0}\,ds\nonumber\\&&{} +\int_0^t\langle
H_\e(z_\e(s)), \phi(s)\rangle_{X^0}ds+ \int_0^t\langle\Sigma dW(s),
\phi(s)\rangle_{X^0}\nonumber
\end{eqnarray}
for $t\in[0,T)$ and $\phi\in C^1(0, T; X_\e^1)$. Moreover if $z^0$
is independent of $W(t)$ with $\mathbf{E}|z^0|^4_{X^0} <\infty $,
then there is positive constant $C_T$, which is independent of $\e$,
such that the following estimates hold:
\begin{equation}\label{est1}
\mathbf{E}|z_\e(t)|^2_{X^0}+\mathbf{E}\int_0^t|z_\e(s)|^2_{X^1_\e}ds\leq
(1+\mathbf{E}|z^0|^2_{X^0})C_T,\;\; {\rm for}  \;\; t\in [0, T]
\end{equation}
and
\begin{equation}\label{est2}
\mathbf{E}\big\{\sup_{t\in[0, T]}|z_\e(t)|^2_{X^0}\big\}\leq
\big(1+\mathbf{E}|z^0|^2_{X^0}+\mathbf{E}\int_0^T|z_\e(s)|^2_{X_\e^1}ds\big)C_T.
\end{equation}
\end{theorem}

\begin{proof}
For any integer $n$, we introduce the following cut-off function
$P^n: \R^+\rightarrow \R^+$ which is a smooth function satisfying
$P^n(x)=1$ if $x<n$ and $P^n(x)=0$ if $x>n+1$. Then nonlinear
function $P^n(|z|)H_\e(z)$ is Lipschitz in both $X^0$ and $X^1$
where $|z|=|u|+\sqrt{\e}|\gamma_1u|$ for $z=(u,
\sqrt{\e}\gamma_1u)$. Now we have the following system with globally
Lipschitz nonlinear term
\begin{equation}
\tilde{z}^n_\e(t)=\mathcal{S}(t)z^0+\int_0^t\mathcal{S}(t-s)
P^n(|\tilde{z}^n_\e(s)|)H_\e(\tilde{z}^n_\e(s))ds+\int_0^t\mathcal{S}(t-s)\Sigma\,dW(s).
\end{equation}
Define a random stopping time $\tau^n(R)$ by
$$
\tau^n(R)=\inf\left\{t>0: |\tilde{z}^n(t)|\geq R \right \}.
$$

Fix arbitrarily a positive number $R<n$ and denote by $\chi_I$ the
characteristic function of the set $I$. Consider the following
integral equation for $t<\tau^n(R)$:
\begin{eqnarray}\label{mild0}
z^n_\e(t)&=&\mathcal{S}(t)z^0+\int_0^{t\wedge\tau^n(R)}
\mathcal{S}(t-s)P^n(|z^n_\e(s)|)H_\e(z^n_\e(s))\,ds\nonumber\\
&&{}+\int_0^{t\wedge\tau^n(R)}\mathcal{S}(t-s)\Sigma\, dW(s)\,.
\end{eqnarray}
Then by the  Theorem 7.4 in \cite{PZ92}, for any $T>0$, the equation
(\ref{mild0}) has a unique solution $z^n_\e\in L^2(\Omega, C(0, T;
X^0)\cap L^2(0, T; X^1))$. Moreover $z_\e^n$ is independent of $n$,
for $n>R$ and satisfies (\ref{SEE}) for $t<\tau^n(R)$.

In the following  we first derive some a priori estimates for
$z_\e$ in $L^2(\Omega, C(0, T; X^0)\cap L^2(0, T; X^1))$.  Then we
  prove the wellposdness of the problem (\ref{SEE}).

Applying the It\^o formula to $|z_\e|^2_{X^0}$ yields
\begin{equation*}
d|z_\e|^2_{X^0}+2\langle\mathcal{A}z_\e, z_\e\rangle_{X^0}dt=2\langle H_\e(z_\e), z_\e\rangle_{X^0}dt+
2\langle\Sigma dW(t), z_\e\rangle_{X^0}+\|\Sigma \|^2_{\mathcal{L}_2^Q}.
\end{equation*}
Here $\|\cdot\|_{\mathcal{L}_2^Q}$ denotes the Hilbert--Schmidt norm
of operator from $Q^{1/2}X^0$ to $X^0$\,. Similar for
$\|\cdot\|_{\mathcal{L}_2^{Q_1}}$ and
$\|\cdot\|_{\mathcal{L}_2^{Q_2}}$ in the following.  By the
assumption ($\mathbf{F}$) we have
\begin{equation*}
\langle H_\e(z_\e), z_\e\rangle_{X^0}=\langle f(u_\e), u_\e\rangle_D
\leq b_1{\rm Mes}(D)
\end{equation*}
where Mes$(D)$ is the Lebesgue measure of the domain $D$.   From
the assumptions ($\mathbf{W}$) and ($\mathbf{\Sigma}$), for any
$\kappa\in (1/2, 1)$, we deduce
\begin{eqnarray*}
\|\Sigma\|^2_{\mathcal{L}^Q_2}
&=&\|\sigma_1\|^2_{\mathcal{L}^{Q_1}_2}+
\|\sigma_2\|^2_{\mathcal{L}^{Q_2}_2}\\
&\leq& \sigma^2_1 trQ_1+\sigma^2_2 tr Q_2\,.
\end{eqnarray*}
  Now combining all the above analysis and (\ref{coercive}) yields
\begin{eqnarray}\label{z2}
d|z_\e|^2_{X^0}+\alpha|z_\e|^2_{X^1}dt\leq \big[\bar{b}+\sigma^2_1
trQ_1+\sigma^2_2 tr Q_2\big]dt+ 2\langle\Sigma\,dW(t),
z_\e\rangle_{X^0}
\end{eqnarray}
where $\bar{b}=2b_1Mes(D)$.
Integrating from $0$ to $t$ and taking expectation on both sides of
the above formula we have the estimate (\ref{est1}) by the Gronwall
inequality \cite{PZ92}. For a further estimate, we apply the It\^o
formula to $|z_\e|^{2m}_{X^0}$, $m>1$, and by calculation similar to
get (\ref{est1}) we have
\begin{equation*}
\frac{d}{dt}\mathbf{E}|z_\e|^{2m}_{X^0}+
\mathbf{E}\big(|z_\e|^{2m-2}_{X^0}|z_\e|^2_{X^1_\e}\big)
\leq C_1\mathbf{E}|z_\e|^{2m-2}_{X^0}+C_2t
\end{equation*}
for some positive constants $C_1$ and $C_2$. Then by the Gronwall
inequality we have
\begin{equation}\label{zm}
\mathbf{E}|z_\e(t)|^{2m}_{X^0}+\mathbf{E}\int_0^t|z_\e(s)|^{2m-2}_{X^0}|z_\e(s)|^2_{X^1_\e}ds\leq
C_T(1+\mathbf{E}|z^0|^{2m}_{X^0}).
\end{equation}

Integrating both sides of (\ref{z2}) from $0$ to $t$, we have
\begin{eqnarray}
\sup_{0\leq t\leq T}|z_\e(t)|^2_{X^0}&\leq &|z^0|^2_{X^0}+
[\sigma^2_1 trQ_1+\sigma^2_2 tr Q_2+\bar{b}]T+\nonumber\\
&&\sup_{0\leq t\leq T}\Big|\int_0^t\langle \Sigma \dot{W}(s),
z_\e(s)\rangle_{X^0}ds\Big|^2+1.\label{zs}
\end{eqnarray}
By the Burkholder-Davis-Gundy inequality \cite{PZ96} and the
assumption ($\mathbf{\Sigma}$) we have
\begin{eqnarray*}
&&\mathbf{E}\sup_{0\leq t\leq T}
\Big|\int_0^t\langle \Sigma \dot{W}(s),  z_\e(s)\rangle_{X^0}ds\Big|^2\\
&\leq&C\mathbf{E}\int_0^T|z_\e(s)|^2_{X^0}\|\Sigma \|^2_{\mathcal{L}_2^Q}ds\\
&\leq&C'\mathbf{E}\int_0^T|z_\e(s)|^2_{X^0}ds
\end{eqnarray*}
for some positive constants $C$ and $C'$. Then by (\ref{est1}) we obtain the estimate (\ref{est2}) from
(\ref{zs}).

Now we continue to prove the wellposedness of the system
(\ref{SEE}). By the assumption of $f$, $z^n_\e$ also satisfies the
estimates (\ref{est1}) and (\ref{est2}) which is independent of
$n$ and $R$. Then we have $\tau(R)\rightarrow \infty$ almost
surely as $R\rightarrow \infty$. For any $T>0$ define $z_\e(t)=
z_\e^n(t)$ for some $n$ and $R$ with $n>R$ and $\tau(R)>T$ almost
surely. Thus $\mathbb{P}(\tau(R)\leq T) = 0$ as $R\to \infty$ for
any $T>0$; for more details about proving global existence, see
\cite{Chow}. The uniqueness and continuity on initial value of
$z_\e$ follows from those of $z_\e^n$. Finally, by the stochastic
Fubini theorem, i.e., Theorem 4.18 in \cite{PZ92}, and the same
discussion of \cite{ChSch04}, we have (\ref{weak1}).

The proof is complete.
\end{proof}

By   Theorem \ref{wellpose} and the definition of $z_\e$ we have
the following corollary.

\begin{coro}\label{est0}
Under the same conditions as in Theorem \ref{wellpose},   for
$t\in [0, T]$, the following estimates hold:
\begin{eqnarray}\label{est3}
&& \mathbf{E}\big(  |u_\e(t)|^2_{L^2(D)}+\e|\gamma_1
u_\e(t)|^2_{L^2(\Gamma_1)}
 \big)+ \nonumber \\
&& \int_0^t \mathbf{E}\big(
|u_\e(t)|^2_{\mathbf{H}_{\Gamma_2}^1(D)}+\e|\gamma_1
u_\e(t)|^2_{\mathbf{B}^1(\Gamma_1)}
 \big)ds \leq
(1+\mathbf{E}|z^0|^2_{X^0})C_T
\end{eqnarray}
and
\begin{eqnarray}\label{est4}
\mathbf{E}\big\{\sup_{t\in[0, T]} [|u_\e(t)|^2_{L^2(D)}+\e|\gamma_1
u_\e(t)|^2_{L^2(\Gamma_1)}] \big\}\leq
(1+\mathbf{E}|z^0|^2_{X^0})C_T\,,
\end{eqnarray}
for some positive constant $C_T$\,, independent of $\e$\,.
\end{coro}

Since the nonlinear term increases polynomially, in order to pass
the limit $\e\rightarrow 0$ in system (\ref{fs1})--(\ref{fs5}),  we
need a priori estimates for $u_\e$ in the space
$\mathbf{H}_{\Gamma_2}^1(D)$\,. In fact we have the following lemma.

\begin{lemma}\label{u1}
Assume that $(\bf{F})$ and
$(\mathbf{\Sigma})$ hold. Let $z^0=(u^0, v^0)\in X^0$ be a
$\big(\mathcal{F}_0, \mathcal{B}(X^0)\big)$-measurable random
variable with $\mathbf{E}|z^0|^2_{X^1}<\infty$. Then for any
$T>0$,
  the solution (\ref{SEE}) $z_\e=(u_\e,
\sqrt{\e}\gamma_1u_\e)\in L^2(\Omega, L^2(0, T; X^1)\cap C(0, T;
X^0))$ and it satisfies the following estimates:
\begin{equation}\label{uv}
\mathbf{E}|u_\e(t)|^2_{\mathbf{H}_{\Gamma_2}^1(D)}+
\mathbf{E}\int_0^t|\D u_\e(s)|^2_{L^2(D)}ds\leq
(1+\mathbf{E}|z^0|^2_{X^1})C_T
\end{equation}
and
\begin{equation}\label{uv1}
\mathbf{E}\sup_{s\in[0, t]}|u_\e(s)|^2_{\mathbf{H}_{\Gamma_2}^1(D)}\leq(1+\mathbf{E}|z^0|^2_{X^1})C_T
\end{equation}
for  any $t\in [0, T]$ and for some positive constant  $C_T$,
independent of $\e$.
\end{lemma}
\begin{proof}
Let $V(z)=a(z,z)$. Applying the It\^o formula to $V(z_\e)$, we have
\begin{eqnarray}\label{u_H1}
&&\frac{1}{2}\frac{d}{dt}V(z_\e)\\
&=&-\langle \mathcal{A}_\e z_\e, \mathcal{A}_\e z_\e\rangle_{X_0}
+\langle H(z_\e), \mathcal{A}_\e z_\e\rangle_{X_0}+\langle
\Sigma \dot{W}, \mathcal{A}_\e z_\e\rangle_{X_0}\nonumber\\ &&+
 \frac{1}{2}\Big[\|\sigma_1\|^2_{\mathcal{L}_2^{A^{\frac{1}{2}}Q_1}}+
 \|\sigma_2\|^2_{\mathcal{L}_2^{Q_2}}\Big]\,.\nonumber
\end{eqnarray}

By definition of $\mathcal{A}_\e$ and (\ref{fs3}) we have
\begin{eqnarray}
\langle H(z_\e)\,, \mathcal{A}_\e z_\e\rangle_{X_0}&=& -\langle
f(u_\e)\,, \D u_\e\rangle_D  \nonumber \\&=&\langle f'(u_\e)\nabla
u_\e\,, \nabla u_\e \rangle_D-
\langle \p_\nu u_\e\,, f(u_\e)\rangle_{\Gamma_1}\nonumber\\
&=&\langle f'(u_\e)\nabla u_\e\,, \nabla u_\e \rangle_D+\langle
\e\dot{u}_\e\,, f(u_\e)\rangle_{\Gamma_1}
+\langle u_\e\,, f(u_\e)\rangle_{\Gamma_1}\nonumber\\
&&-\langle \sqrt{\e}\sigma_2\dot{W}_2\,, f(u_\e)
\rangle_{\Gamma_1}\,.\label{fD}
\end{eqnarray}

Also by the It\^o formula
$$
\langle \e\dot{u}_\e\,,
f(u_\e)\rangle_{\Gamma_1}=\e\frac{d}{dt}\langle  F(u_\e)\,,
1\rangle_{\Gamma_1} -\frac{\sigma^2_2}{2}tr \big( D_{uu}\langle
F(u_\e)\,, 1\rangle_{\Gamma_1}
Q_2^{\frac{1}{2}}(Q_2^{\frac{1}{2}})^*\big)\,.
$$

Then we can have from (\ref{u_H1})--(\ref{fD})
\begin{eqnarray}\label{e}
&&\frac{1}{2}\frac{d}{dt}\Big [V(z_\e)-
2\e\langle F(u_\e)\,, 1 \rangle_{\Gamma_1} \Big]\\
&=&-|\mathcal{A}_\e z_\e|^2_{X_0}-\langle \Sigma\dot{W}\,,
\mathcal{A}_\e z_\e\rangle_{X_0}+
\frac{1}{2}\Big[\|\sigma_1\|^2_{\mathcal{L}_2^{A^{\frac{1}{2}}Q_1}}+
 \|\sigma_2\|^2_{\mathcal{L}_2^{Q_2}}\Big]\nonumber\\
&&+\langle f'(u_\e)\nabla u_\e\,, \nabla u_\e  \rangle_D+\langle
u_\e\,, f(u_\e)\rangle_{\Gamma_1}-\frac{\sigma^2_2}{2}tr \big(
D_{uu}\langle F(u_\e)\,, 1\rangle_{\Gamma_1}
 Q_2^{\frac{1}{2}}(Q_2^{\frac{1}{2}})^*\big).\nonumber
\end{eqnarray}
By assumptions ($\mathbf{F}$) and ($\mathbf{\Sigma}$), taking
expectation on both  sides of (\ref{e}), applying the Cauchy
inequality and noticing that
\begin{equation}\label{ineq}
|u_\e|^2_{L^2(\Gamma_1)}\leq
C_\varepsilon|u_\e|^2_{L^2(D)}+\varepsilon|u_\e|^2_{\mathbf{H}^1_{\Gamma_2}(D)}
\end{equation}
for $\varepsilon>0$, we can have (\ref{uv}) by the Gronwall
inequality and Corollary \ref{est0} by taking $\varepsilon>0$
small enough.


Integrating both sides of (\ref{e}) from $0$ to $t$, and using the
Cauchy inequality, we conclude that
\begin{eqnarray*}
&&\frac{1}{2}\mathbf{E}\sup_{t\in [0, T]}|u_\e(t)|^2_{\mathbf{H}_{\Gamma_2}^1(D)}\\
&\leq& \frac{1}{2}\mathbf{E}|u^0|^2_{\mathbf{H}^1_{\Gamma_2}(D)}+ C\mathbf{E}\Big[\int_0^T|u_\e(s)|^2_{L^2(D)}ds+
\int_0^T|u_\e(s)|^2_{\mathbf{H}^1_{\Gamma_2}(D)}ds+1\Big]
\end{eqnarray*}
for some positive constant $C$ depending only on $\tilde{b}$, $Q_1$,
$Q_2$, $\sigma_1$ and $\sigma_2$. Then we   have (\ref{uv1}) by
Corollary \ref{est0}\,.

The proof is complete.
\end{proof}

\vskip 0.8cm


\section{Tightness of the distributions of solutions}

We intent  to investigate the limit of the solution $u_\e$ of
stochastic system (\ref{fs1})--(\ref{fs5}) as $\e\rightarrow 0$ in
the sense of distribution. For this purpose, in this section, we
establish results on tightness of the distributions of the
solutions. Let $\mu_\e$ be the distribution of $u_\e$, which
generates a Radon probability measure on the following metric space
$$
\mathcal{H} := L^2(0, T; L^2(D)) \cap C(0, T; H^{-1}(D)).
$$
 Now we prove that the family of distributions   $\{\mu_\e \}$ is tight in the
space $\mathcal{H}$.

We apply the a priori estimates in the  preceding section to obtain
the tightness of $\{ \mu_\e\}$.

 First by the property of Wiener process, for some $\rho\in (0, \frac{1}{2})$
\begin{equation}\label{Me}
\mathbf{E}\sup_{|t-r|\leq \tau}\frac{|\sigma_1 W(t)-\sigma_1W(r)|_{L^2(D)}
}{|t-r|^\rho}\leq \sigma_1C_T.
\end{equation}
Then, by estimate (\ref{est1}) and (\ref{uv1}), for any given $\delta>0$,
there is a positive constant $C^\delta_T$
 such that
 \begin{equation*}
 \mathbf{P}\{ A_\delta \}>1-\delta
 \end{equation*}
with
\begin{align}
A_\delta=\big\{\omega\in \Omega:\sup_{0\leq t\leq T}|u_\e(t)|_{\mathbf{H}^1_{\Gamma_2}(D)} \leq C^\delta_T&,&
\sup_{|t-r|\leq \tau}\frac{|\sigma_1 W(t)-\sigma_1W(r)|_{L^2(D)}
}{|t-r|^\rho}\leq C_T^\delta \nonumber \\ {\rm and }\;
\int_0^T|u_\e(t)|^2_{H^1_{\Gamma_2}(D)}dt\leq C_T^\delta \big\}.\label{Adelta}
\end{align}
For any $\varphi\in H_0^1(D)$,  by (\ref{fs1})--(\ref{fs5}), we have
\begin{equation*}
d\langle u_\e(s), \varphi \rangle_D =a(u_\e, \varphi)ds +\langle
f(u_\e), \varphi\rangle_Dds+ \langle \sigma_1d{W}_1(s),
\varphi \rangle_D.
\end{equation*}
By the definition of $a(z,\bar{z})$, we have
\begin{equation}\label{a1}
 \sup_{|t-r|\leq \tau}\left|\int_r^ta(u_\e, \varphi)ds\right|\leq
M\Big[\int_0^T|u_\e(s)|^2_{\mathbf{H}^1_{\Gamma_2}(D)}ds\Big]^{1/2}|\varphi|_{H^1_0(D)}\sqrt{\tau}
\end{equation}
where $M$ is defined in (\ref{e:a})\,. By the assumption
$\mathbf{(F)}$ and the embedding $\mathbf{H}_0^1(D)$ into $L^6(D)$
we have
\begin{equation}\label{f1}
 \sup_{|t-r|\leq \tau}\left|\int_r^t\langle f(u_\e), \varphi
\rangle_Dds\right|\leq
C'\Big[1+\int_0^T|u_\e(s)|^6_{\mathbf{H}^1_{\Gamma_2}(D)}ds
\Big]^{1/2} |\varphi|_{H^1_0(D)}\sqrt{\tau}
\end{equation}
for some positive constant $C'$\,. Then if $\omega\in A_\delta$,  by
the definition of $A_\delta$ and (\ref{a1})--(\ref{f1})
$$
|u_\e|_{C^{\rho}(0,T; H^{-1}(D) )\cap L^2(0, T;
\mathbf{H}^1_{\Gamma_2}(D))}\leq C(T, \delta)
$$
for some positive constant $C(T, \delta)$\,.  Define set
$K_\delta\subset\mathcal{H}$ as
$$
K_\delta=\big\{u\in \mathcal{H}:|u_\e|_{C^{\rho}(0,T; H^{-1}(D)
)\cap L^2(0, T; \mathbf{H}^1_{\Gamma_2}(D))}\leq C(T, \delta)
\big\}.
$$
Then by the compact embedding of $C^\rho(0, T;H^{-1}(D) )\cap L^2(0, T; \mathbf{H}^1_{\Gamma_2}(D))$
into $\mathcal{H}$, $K_\delta$ is compact in $\mathcal{H}$. And by the definition of $A_\delta$ and above
analysis we have
\begin{equation*}
\mathbf{P}\{ u_\e\in K_\delta \}>1-\delta.
\end{equation*}
Then we proved the following result.
\begin{theorem} \label{tightness} (\textbf{Tightness of distributions of
solutions})\\
The family of distributions of the solutions, $\{ \mu_\e \}$, is
tight in the space $\mathcal{H}$.
\end{theorem}


\vskip 0.8cm


\section{Effective   dynamics}\label{s5}

In this section we pass the limit of $\e\rightarrow 0$ in
(\ref{weak}) and obtain the limiting system in a certain sense to be
specified below. We always assume that
$\mathbf{E}|z^0|^2_{X^1}<\infty$ and
$\mathbf{E}|z^0|^4_{X^0}<\infty$ in the following.


Since $\{ \mu_\e \}$ is tight in the space $\mathcal{H}$ (defined in
the beginning of the last section), for any $\delta>0$ there is a
compact set $K_\delta\subset \mathcal{H}$ such that
$$
\mathbf{P}\{u_\e\in K_\delta\}>1-\delta.
$$
Here $K_\delta$ is chosen as a family of decreasing sets with
respect to $\delta$, i.e. $K_\delta\subset K_{\delta'}$ for any
$\delta\geq\delta'>0$. Moreover by the analysis of last section we can choose the
set $K_\delta$ with some positive constant $C_T^\delta$, depending
on $T$ and $\delta$, such that
 $\sup_{t\in [0, T]}|u_\e(t)|^2_{\mathbf{H}^1_{\Gamma_2}(D)}\leq C^\delta_T$  for $u_\e\in
K_\delta$.

Then Prohorov theorem and Skorohod embedding theorem (\cite{PZ92})
assure  that for any sequence $\{\e_ j\}_j$ with $\e_j\rightarrow 0$
as $j\rightarrow \infty$, there exist subsequence $\{\e_ {j(k)}\}$,
random elements  $\{u^*_{\e_{j(k)}}\}\subset \mathcal{H}$, $u^*\in
\mathcal{H}$, $u^{0*}\in L^2(D)$ and $L^2(D)$-valued Wiener process
$W^*_1$, $L^2(\Gamma_1)$-valued Wiener process $W^*_2$ defined on a
new probability space $(\Omega^*, \mathcal{F}^*, \mathbb{P}^*)$,
such that
$$
\mathcal{L}(u^*_{\e_ {j(k)}})=\mathcal{L}(u_{\e_{j(k)}})
$$
and
$$
u^*_{\e_ {j(k)}}\rightarrow u^*\;\;in\;\;\mathcal{H}\;\; as
\;\;k\rightarrow \infty,
$$
for almost all $\omega\in \Omega^*$. Moreover $u^*_{\e_ {j(k)}}$
solves system (\ref{fs1})--(\ref{fs5}) with $W_1$ and $W_2$ are
replaced by Wiener process $W^*_{1k}$ and $W^*_{2k}$ respectively
defined on probability space $(\Omega^*, \mathcal{F}^*,
\mathbb{P}^*)$ with same distribution as $W_1$ and $W_2$ for any
$k$. And $u^0$ is replaced by a random variable $u_k^{0*}$ with
$\mathcal{L}(u^0)=\mathcal{L}(u_k^{0*})$. And for $\mathbb{P}^*$-
almost all $\omega\in\Omega^*$,
$|u_k^{0*}-u^{0*}|_{L^2(D)}\rightarrow 0$ and
$$
\sup_{0\leq t \leq T}|W^*_1(t)-W^*_{1k}(t)|_{L^2(D)}\rightarrow
0,\;\;\sup_{0\leq t \leq
T}|W^*_2(t)-W^*_{2k}(t)|_{L^2(\Gamma_1)}\rightarrow 0
$$
for $k\rightarrow \infty$. Now we will determine the limiting
equation satisfied by $u^*$ and the limiting equation is independent
of $\e$. In fact we will prove that  $u^*$ solves
(\ref{avsys1})--(\ref{avsys4}) with $u^0$ and $W_1$ are replaced by
$u^{0*}$ and $W_1^*$ respectively.

We will pass the limit $\e\rightarrow 0$ in (\ref{weak}) for
$u^*_{\e_{j(k)}}$. For the nonlinear term $f(u)$ increases
polynomially, in order to pass the limit in $f(u^*_{\e_{j(k)}})$ we
restrict $u_\e$ in a bounded set in $\mathbf{H}^1_{\Gamma_1}(D)$.
However this is impossible for stochastic process $u^*_{\e_{j(k)}}$
which converges just in space $\mathcal{H}$.  For this define a new
probability space $(\Omega^*_\delta, \mathcal{F}^*_\delta,
\mathbb{P}^*_\delta)$ as
\begin{equation*}
 \Omega^*_\delta=\{\omega'\in \Omega: u^*_{\e_{j(k)}}\in
K_\delta\},
\end{equation*}
\begin{equation*}
\mathcal{F}^*_\delta=\{S\cap \Omega^*_\delta: S\in\mathcal{F}^*\}
\end{equation*}
and
$$
\mathbb{P}^*_\delta(S)=\frac{\mathbb{P}^*(S\cap\Omega^*_\delta)}{\mathbb{P}(\Omega^*_\delta)},
\;\;{\rm for} \;\; S\in\mathcal{F}^*_\delta.
$$
Denote by $\mathbb{E}^*_\delta$ the expectation operator with
respect to $\mathbb{P}^*_\delta$.   It is clear that
$\mathbf{P}(\Omega^*\setminus\Omega^*_\delta)\leq\delta$.

Since the distribution of  $u^*_{\e_{j(k)}}$ is same as that of
$u_{\e_{j(k)}}$\,,  $u^*_{\e_{j(k)}}$  converges to $u^*$ in  space
$L^2(0, T; L^2(D))\cap C(0, T; \mathbf{H}^{-1}(D))$ under the usual
metric for $\omega\in\Omega^*_\delta$. Here the usual metric in the
intersection $X\cap Y$, of two metric spaces $(X, d_X)$ and $(Y,
d_Y)$, is the metric $d:= d_X+d_Y$ or equivalently $d:= \max\{d_X,
d_Y \}$.


Now in the following we determine the limiting system satisfied by
$u^*$ restricted on probability space $( \Omega^*_\delta,
\mathcal{F}^*_\delta, \mathbb{P}^*_\delta)$. By the weak solution in
the sense of (\ref{weak1}), for any $\psi\in C^1(0,T;C^\infty(D))$
with $\psi(T)=0$ we have
\begin{eqnarray}\label{pslimit}
 &&-\int_0^T\langle u^*_{\e_{j(k)}}(t), \dot{\psi} \rangle_Ddt-
 \e  \int_0^T\langle \gamma_1 u^*_{\e_{j(k)}}(t),\dot{\psi} \rangle_{\Gamma_1} dt \\
&=&\int_0^T a(u^*_{\e_{j(k)}}, \psi)dt +\int_0^T\langle
f(u^*_{\e_{j(k)}}), \psi\rangle_Ddt+ \int_0^T\langle
\sigma_1\dot{W}^*_{1k}(t),
\psi \rangle_Ddt+ \nonumber\\
 && \sqrt{\e}\int_0^T\langle \sigma_2\dot{W}^*_{2k}(t),
 \psi\rangle_{\Gamma_1} dt+\langle u_k^{0*}, \psi(0)\rangle_D+
 \e\langle \gamma_1 u_k^{0*}, \gamma_1 \psi(0)
 \rangle_{\Gamma_1},  \nonumber
\end{eqnarray}
on  $( \Omega^*_\delta, \mathcal{F}^*_\delta, \mathbb{P}^*_\delta) $
.

 We consider the terms in (\ref{pslimit}) respectively. Since
$u^*_{\e_{j(k)}}$ converges weakly  to $u^*$ in $L^2(0, T;
\mathbf{H}^1_{\Gamma_2}(D))$
\begin{equation}\label{1}
\int_0^T a(u^*_{\e_{j(k)}}, \psi)dt\rightarrow \int_0^T a(u^*,
\psi)dt, \;\; \e\rightarrow 0
\end{equation}
for all $\omega\in\Omega_\delta^*$. By the definition of $
\Omega^*_\delta$ and the choice of $K_\delta$,
$|u^*_{\e_{j(k)}}(t)|_{\mathbf{H}^1_{\Gamma_2}(D)}$ is bounded
uniformly in $t\in [0, T]$ and $\e\in (0, 1]$. By the embedding of
$\mathbf{H}^1_{\Gamma_2}(D)$ into $L^{\frac{2N}{N-2}}(D)$,
$|f(u^*_{\e_{j(k)}})|_{L^2(0,\; T; L^2(D))}$ is bounded.  Then by
Lemma \ref{Lions} and assumption ($\mathbf{F}$),
$f(u^*_{\e_{j(k)}})$ converges weakly to $f(u^*)$ in $L^2(0, T;
L^2(D))$ which means for $\omega\in\Omega^*_\delta$
\begin{equation}\label{2}
\int_0^T\langle f(u^*_{\e_{j(k)}}), \psi\rangle_Ddt\rightarrow
\int_0^T\langle f(u^*), \psi\rangle_Ddt, \;\; \e\rightarrow 0.
\end{equation}

By assumption ($\mathbf{\Sigma}$) and  the property of stochastic
integral, see also Lemma 3.1 of \cite{GyKry96}
\begin{equation}\label{}
\int_0^T\langle \sigma_1\dot{W}^*_{1k}(s),
\psi(s) \rangle_Dds\rightarrow \int_0^T\langle
\sigma_1\dot{W}^*_{1}(s), \psi(s) \rangle_Dds
\end{equation}
in $\mathbb{P}^*_\delta$-probability  and
\begin{equation}\label{}
\sqrt{\e}\int_0^T\langle
\sigma_2\dot{W}^*_{2k}(s), \psi(s)
\rangle_{\Gamma_1}ds\rightarrow 0
\end{equation}
in $\mathbb{P}^*_\delta$-probability.
 Moreover for $\omega\in\Omega^*_\delta$
\begin{eqnarray}\label{5.5}
&&\big|\e  \int_0^T\langle \gamma_1 u^*_{\e_{j(k)}}(t),\dot{\psi}
\rangle_{\Gamma_1}dt\big|^2\nonumber \\
&\leq & \e^2 T|\psi|_{C^1(0, T;
 H^1(D))} \int_0^T|u^*_{\e_{j(k)}}|^2_{\mathbf{H}^1_{\Gamma_2}(D)}dt \nonumber \\
 &\rightarrow& 0,\;\;\e\rightarrow 0. \label{5}
\end{eqnarray}
Then combining the above analysis in (\ref{1})--(\ref{5}) and by the
density argument we could conclude
\begin{eqnarray}\label{e0}
 &&-\int_0^T\langle u^*(t), \dot{\psi} \rangle_Ddt \\
&=&\int_0^T a(u^*, \psi)dt +\int_0^T\langle f(u^*), \psi\rangle_Ddt+
\int_0^T\langle \sigma_1\dot{W}^*_{1}(t),
\psi \rangle_Ddt+ \nonumber\\
 && +\langle u^{0*}, \psi(0)\rangle_D \nonumber.
\end{eqnarray}

Integrating by parts in (\ref{e0}) we see  that the limiting
function $u^*$ satisfies the following system with deterministic
boundary condition
\begin{eqnarray}\label{sys1}
du^*&=&\big [\D u^* +f(u^*) \big]dt+
\sigma_1\,dW^*_{1}(t)\;\;in\;\; D,\\
\p_\nu u^*+u^*&=&0\;\; on \;\; \Gamma_1, \label{sys2}\\
u^*&=&0\;\; on\;\; \Gamma_2,\\
u^*(0)&=&u^{0*}\;\; in\;\; D, \label{sys4}
\end{eqnarray}
on the probability space $( \Omega^*_\delta, \mathcal{F}^*_\delta,
\mathbb{P}^*_\delta)$, which  has static boundary condition. For the
above system (\ref{sys1})--(\ref{sys4}) we can rewrite in the
following abstract form
\begin{equation}\label{abssys}
du^*=[-A_0u^*+f(u^*)]dt+\sigma_1\,dW^*_{1}(t),\;\;
 u^*(0)=u^{0*}
\end{equation}
where the  operator $-A_0$ is the Laplace operator with the Robin
boundary condition. The corresponding bilinear form is
$$
a_0(v_1, v_2)=\int_D\nabla v_1 \nabla v_2\,
dx+\int_{\Gamma_1}(\gamma_1v_1)(\gamma_1v_2)\,d\Gamma_1.
$$
Then by Theorem 7.4 of \cite{PZ92}, and a similar analysis in the
proof of Theorem \ref{wellpose}, for any $T>0$ system
(\ref{efsys1})--(\ref{efsys4}) has a unique solution $u^*\in
L^2\big( \Omega^*_\delta, L^2(0, T; \mathbf{H}^1_{\Gamma_2}(D))\cap
C(0, T; L^2(D))\big)$ in the sense of (\ref{e0}).

Then by the arbitrariness of the choice of $\delta$,
\begin{equation}\label{wstar}
u^*_{\e_{j(k)}} \;\; {\rm  converges\;\; in \;\;\mathbb{P}^*-
probability\;\; to \;\;} u^*
\end{equation}
which solves
\begin{eqnarray}\label{efsys1}
du^*&=&\big [\D u^* +f(u^*) \big]dt+
\sigma_1 \; dW^*_1(t)\;\;in\;\; D,\\
\p_\nu u^*+u^*&=&0\;\; on \;\; \Gamma_1, \label{efsys2}\\
u^*&=&0\;\; on\;\; \Gamma_2,\\
u^*(0)&=&u^{0*}\;\; in\;\; D, \label{efsys4}
\end{eqnarray}
on the probability space  $( \Omega^* , \mathcal{F}^* , \mathbb{P}^*
)$. In order to obtain the convergence in probability of $u_\e$, we
need the following lemma from \cite{GyKry96}.
\begin{lemma}\label{convprob}
Let $Z_n$ be a sequence of random elements in a Polish space
$(\mathbb{X}, d)$ equipped with Borel $\sigma$-algebra. Then $Z_n$
converges in probability to an $\mathbb{X}$-valued random element if
and only if for every pair subsequences $Z_l$ and $Z_m$, there
exists a subsequence $v_k:=(Z_{j(k)}, Z_{m(k)})$ converging weakly
to a random element $v$ supported on the diagonal $\{(x,
y)\in\mathbb{X}\times\mathbb{X}: x=y \}$.

\end{lemma}

Here we call an $\mathbb{X}$-valued random variable $X_n$ converges
weakly to $X$ if
\begin{equation*}
\mathbf{E}f(X_n)=\int_\mathbb{X}f(x)\mathbb{P}_n(dx)\rightarrow
\int_\mathbb{X}f(x)P(dx)=\mathbf{E}f(X)
\end{equation*}
with $\mathbb{P}_n=\mathcal{L}(X_n)$ and
$\mathbb{P}=\mathcal{L}(X)$. Notice that convergence in probability
implies weakly convergence, see \cite{Bill, Dudley}. Then by the
uniqueness property of solution for equations
(\ref{avsys1})--(\ref{avsys4}) which can be proved similarly by that
for equations (\ref{sys1})--(\ref{sys4}), we can formulate the main
result in this section by Lemma \ref{convprob}.
\begin{theorem}\label{main1} (\textbf{Effective   system})\\
Assume that  conditions  $(\mathbf{F})$ and $(\mathbf{\Sigma})$ are
satisfied. Let $u^0\in L^2(D)$ be a $(\mathcal{F}_0,
\mathcal{B}(L^2(D)))$- measurable random variable, which is
independent of $W(t)$, with
$\mathbf{E}|u^0|^2_{\mathbf{H}^1_{\Gamma_2}(D)}<\infty$ and
$\mathbf{E}\big[|u^0|^4_{L^2(D)}+ |\gamma_1 u^0|^4_{L^2(\Gamma_1)}
\big]<\infty$. Then for any $T>0$, the solution $u_\e$ of the
stochastic system (\ref{fs1})--(\ref{fs5}) converges to $u$, which
is the solution of the effective limiting system
(\ref{avsys1})--(\ref{avsys4})\,, in probability in
 space $\mathcal{H}$:
$$
\mathcal{H}=L^2(0, T; L^2(D)) \cap C(0, T; H^{-1}(D)).
$$
\end{theorem}


\vskip 0.6cm


\section{Normal deviations}\label{s6}
We have proved $u_\e$ approaches $u$ in probability in space
$\mathcal{H}$, namely, the difference $u_\e-u$ tends to $0$ in
probability in space $\mathcal{H}$ as $\e\rightarrow 0 $. In this
section we consider the order of $u_\e-u$ in $\e$ as $\e \rightarrow
0$, that is, the normal deviations of $u_\e$ away from the effective
solution $u$. We prove that the order is $\frac{1}{2}$ and the
normalized difference $\frac{1}{\sqrt{\e}}(u_\e-u)$ converges in an
appropriate function space. In this section we   assume the assumption in Theorem \ref{main1}.

Denote by $v_\e=\frac{1}{\sqrt{\e}}(u_\e-u)$. Then we have the
following initial boundary value problem for $v_\e$
\begin{eqnarray}
\dot{v}_\e&=&\D v_\e+\frac{1}{\sqrt{\e}}\big(f(u_\e)-f(u)\big),\;\; in\;\; D \label{ve1}\\
\p_\nu v_\e+v_\e&=&-\sqrt{\e}\dot{u}_\e+\sigma_2 \dot{W}_2,\;\; on
\;\;\Gamma_1 \label{ve2}\\
v_\e&=&0,\;\; on \;\;\Gamma_2 \label{ve3} \\
v_\e(0)&=&0. \label{ve4}
\end{eqnarray}

As $\e\rightarrow 0 $ we expect $v_\e$ converges in some sense to
the solution $v$ of the following linear system
\begin{eqnarray}
dv&=&[\D v+f'(u)v]dt,\;\; in\;\; D \label{v1}\\
\p_\nu v+v&=&\sigma_2 \dot{W}_2,\;\; on\;\;
\Gamma_1 \label{v2}\\
v&=&0,\;\; on \;\;\Gamma_2 \label{v3} \\
v(0)&=&0. \label{v4}
\end{eqnarray}
Note that the limiting system (\ref{v1})--(\ref{v4}) contains a
static  boundary with random force (but not dynamical). For the
wellposedness of the above two systems we follow the approach of
\cite{DaPrato3, PZ96}; see also \cite{Mas95}. Let $\mathcal{N}$ be a
linear bounded operator from $L^2(\Gamma_1)$ into $L^2(D)$ defined
as the solution of following problem
\begin{equation}\label{Neum}
ry-\D y=0\;\; in\;\; D, \;\; \p_\nu y+y = g\;\; on \;\; \Gamma_1
\end{equation}
with $r\in\R$ such that (\ref{Neum}) has a unique solution
$y=\mathcal{N} g$ for any $g\in L^2(\Gamma_1)$. Here $\mathcal{N}$
is called the Nenumann mapping. For our problem let
$g=g_\e=-\sqrt{\e}\dot{u}_\e+\sigma_2 \dot{W}_2$ and $g_0=\sigma_2
\dot{W}_2$. However it is easy to see that $g_\e$ and $g_0$ are not
in the space $L^2(\Gamma_1)$. Fortunately, we can extend
$\mathcal{N}$ to a bounded linear operator from
$H^{\varrho}(\Gamma_1)$ to $D^\varepsilon_A$ which is the domain of
the operator $(rI-A)^\varepsilon$ with
$0<\varepsilon<\frac{\varrho}{2}+\frac{3}{4}$,
$\varrho>-\frac{3}{2}$, see \cite{Lion1} or example 3.1 in
\cite{Mas95}. Here $A$ is a second order differential operator
defined on $\{u\in H^2(D), \p_\nu u=0\}$ with $Au=\D u$. Denote by
$S(t)$ the $C_0$ semigroup generated by the linear
 operator $A$. Then we can write the solution of
(\ref{ve1})--(\ref{ve4}) and (\ref{v1})--(\ref{v4}) respectively in
the following mild sense

\begin{eqnarray}
v_\e(t)&=&-\int_0^tAS(t-r)\mathcal{N}(\gamma_1v_\e)(r)dr+\sigma_2\int_0^tAS(t-r)\mathcal{N}dW_2(r)+\nonumber\\
&&\frac{1}{\sqrt{\e}}\int_0^tS(t-r)\big[f(u_\e(r))-f(u(r))\big]dr+\nonumber\\
&&\sqrt{\e}\int_0^tAS(t-r)\mathcal{N}(\dot{\gamma_1u}_\e)(r)dr\label{v-mild}
\end{eqnarray}
and
\begin{eqnarray*}
v(t)&=&-\int_0^tAS(t-r)\mathcal{N}(\gamma_1v)(r)dr+\sigma_2\int_0^tAS(t-r)\mathcal{N}dW_2(r)+\\
&&\int_0^tS(t-r)f'(u(r))v(r)dr.
\end{eqnarray*}

By the Example 3.1 of \cite {Mas95} for any $T>0$, there exist
functions $v_\e$ and $v$, both in $C(0, T; L^2(\mathbf{\Omega},
L^2(D)))$, which are unique mild solutions of
(\ref{ve1})--(\ref{ve4}) and (\ref{v1})--(\ref{v4}), respectively.

\begin{remark}
For a special one-dimensional case on the domain $D=(-1,1)$,
 the solution is proved earlier \cite{FW92} in a special weighted space
$\hat{C}_a\subset C_a([0, T]\times (-1, 1))$\,. Here $C_a([0,
T]\times (-1, 1))$ consisting of  all continuous functions
$u(t,x)$\,, $t\in [0, T]$\,, $|x|<1$\,, such that
$\lim_{x\rightarrow\pm1}a(x)u(t,x)=0$ uniformly in $t\in[0, T]$\,.
Here weighted function $a(x)$ can be chosen as $(1-x^2)^\alpha$\,,
$0<\alpha<1$\,. Then $\hat{C}_a\subset C_a([0, T]\times (-1, 1))$
consisting of $u(t, x)$ such that $h(t, x) = \int_0^tu(s, x)\,ds$
has uniform limit in $t\in[0, T]$ when $x\rightarrow 1$ and when
$x\rightarrow-1$\,.
\\
\end{remark}

Let $\nu_\e$ be the distributions of $v_\e$ in the space $L^2(0, T; L^2(D))$.
For our purpose in the following, we prove the tightness of $\nu_\e$. First we should
derive a further a priori estimate  for $v_\e$. As pointed out in \cite{Mas95} the It\^o
formula cannot be used for the Lyapunov function $V(x)=|x|^p$,
$p>0$. We treat $v_\e$ in the mild sense (\ref{v-mild}).
For any $\delta>0$, we still consider $\omega\in A_\delta$ which defined
by (\ref{Adelta}).

\begin{lemma}\label{H-lambda}
Let $u^0\in L^2(D)$ be a $(\mathcal{F}_0, \mathcal{B}(L^2(D)))$-
measurable random variable, which is independent of $W(t)$, with
$\mathbf{E}|u^0|^{2}_{\mathbf{H}^1_{\Gamma_2}(D)}<\infty$ and
$\mathbf{E}\big[|u^0|^4_{L^2(D)}+ |\gamma_1
u^0|^4_{L^2(\Gamma_1)} \big]<\infty$. Then for any $T>0$, there
exist a $\lambda>0$ and    a positive constant $C_T$ such that
 \begin{equation*}
 \mathbf{E}[\chi_{A_\delta}|v_\e|^2_{L^2(0, T; \mathbf{H}_{\Gamma_2}^\lambda(D))} ]\leq C_T.
 \end{equation*}
Here $\chi_{A_\delta}(\omega)=1$ for $\omega\in A_\delta$ and $\chi_{A_\delta}(\omega)=0$
otherwise.

\end{lemma}
\begin{proof}
By the similar estimates in the proof for Proposition 2.2 of
\cite{Mas95} and estimates (\ref{zm}) with $m=p$\,, noticing that
the initial value is zero,  we have
 \begin{equation*}
 \mathbf{E}[\chi_{A_\delta}|v_\e|^q_{L^q(0, T; \mathbf{H}^\lambda_{\Gamma_2}(D))}]\leq C_T
 \end{equation*}
for some $q>2$. Then by the H\"older inequality we have the result.
\end{proof}

%

Furthermore multiplying $\varphi\in C_0^\infty(D)$ to both sides
of (\ref{ve1}) yields
\begin{eqnarray*}
&&\int_0^T\big\langle \dot{v}_\e(t),
\varphi\big\rangle_Ddt\\&=&-\int_0^T\big\langle \nabla v_\e(t),
\nabla
\varphi\big\rangle_Ddt+\frac{1}{\sqrt{\e}}\int_0^T\big\langle(f(u_\e(t))-f(u(t))),
\varphi \big\rangle_Ddt.
\end{eqnarray*}
Then by Lemma \ref{H-lambda} we deduce
\begin{equation}\label{ven1}
\mathbf{E}[\chi_{A_\delta}|\dot{v}_\e(t)|^2_{L^2(0, T;
\mathbf{H}^{\lambda-1}(D))}]\leq C_T
\end{equation}
for some positive constant $C_T$. Then by the Chebyshev
inequality and Lemma \ref{comp}, for any $\delta>0$ there is a
compact subset $N_\delta\subset L^2(0, T; L^2(D))$ such that
$$
\mathbf{P}\{v_\e\in N_\delta \}>1-\delta.
$$
That is, the probability measure sequence $\{\nu_\e\}$ is tight in
space $L^2(0, T; L^2(D))$.  Then Prohorov theorem and Skorohod
embedding theorem (\cite{PZ92}) assure  that for any sequence
$\{\e_ j\}_j$ with $\e_j\rightarrow 0$ as $j\rightarrow \infty$,
there exist subsequence $\{\e_ {j(k)}\}$, random elements
$\{\bar{v}_{\e_{j(k)}}\}\subset L^2(0, T; L^2(D))$, $\bar{v}\in
L^2(0, T; L^2(D))$, $L^2(\Gamma_1)$-valued Wiener process
$\overline{W}_{2k}$ defined on a new probability space
$(\overline{\Omega}, \overline{\mathcal{F}},
\overline{\mathbb{P}})$, such that
$$
\mathcal{L}(\bar{v}_{\e_ {j(k)}})=\mathcal{L}(v_{\e_{j(k)}})
$$
and
$$
\bar{v}_{\e_ {j(k)}}\rightarrow \bar{v}\;\;in\;\;L^2(0, T;
L^2(D))\;\; as \;\;k\rightarrow \infty,
$$
for almost all $\omega\in \overline{\Omega}$. $\bar{v}_{\e_{j(k)}}$
solves (\ref{ve1})--(\ref{ve4}) with $W_2$ replaced by
$\overline{W}_{2k}$. And for almost all $\omega\in\overline{\Omega}$
$$
\sup_{0\leq t\leq
T}|\overline{W}_{2k}-\overline{W}_2|_{L^2(D)}\rightarrow 0,\;\;
k\rightarrow \infty.
$$
Moreover
$\bar{v}_{\e_{j(k)}}=\frac{1}{\sqrt{\e_{j(k)}}}(\bar{u}_{\e_{j(k)}}-\bar{u})$
for some random elements $\bar{u}_{\e_{j(k)}}$, $\bar{u}\in
\mathcal{H}$ with
$\mathcal{L}(\bar{u}_{\e_{j(k)}})=\mathcal{L}(u_{\e_{j(k)}})$ and
$\mathcal{L}(\bar{u})=\mathcal{L}(u)$.


In order to pass limit $\e\rightarrow 0$ in $f'$, by the same
approach of \S \ref{s5} we define the following new probability
space $(\overline{ \Omega }_\delta, \overline{\mathcal{F}}_\delta,
 \overline{\mathbb{P}}_\delta)$ as
\begin{equation*}
\overline{\Omega}_\delta=\{\omega'\in\overline{\Omega}:
v_\e(\omega')\in N_\delta\},
\end{equation*}
\begin{equation*}
\overline{\mathcal{F}}_\delta=\{S\cap\overline{\Omega}_\delta:
S\in\mathcal{F}_\delta\}
\end{equation*}
and
$$
\overline{\mathbb{P}}_\delta(S)=
\frac{\mathbb{P}_\delta(S\cap\overline{\Omega}_\delta)}{\mathbb{P}_\delta(\overline{
\Omega}_\delta)}, \;\;{\rm for} \;\; S\in\mathcal{F}_\delta.
$$
  Now we restrict the
system (\ref{ve1})--(\ref{ve4}) on the probability space
$(\overline{ \Omega}_\delta, \overline{\mathcal{F}}_\delta,
 \overline{\mathbb{P}}_\delta)$.   By the definition of $\overline{\Omega}_\delta$
and the discussion in \S \ref{s5} for almost all
$\omega\in\overline{\Omega}_\delta$
 \begin{equation}\label{as}
 \bar{u}_{\e_{j(k)}}(t,x)\rightarrow \bar{u}(t,x), \;\; {\rm almost\;\; everywhere\;\;in
 } \;\; [0, T]\times D.
 \end{equation}
And for $\omega\in \overline{ \Omega }_\delta$,
$\bar{v}_{\e_{j(k)}}$ converges to $\bar{v}$ almost surely on $[0,
T]\times D$,
$$
\bar{v}_{\e_{j(k)}}\rightarrow \bar{v}\;\; weakly\;\; in\;\; L^2(0,
T; \mathbf{H}^\lambda_{\Gamma_2}(D)),\; as \;k\rightarrow \infty
$$
and
 $$
\bar{v}_{\e_{j(k)}}\rightarrow \bar{v}\;\;strongly \;\;in \;\;
L^2(0, T; L^2(D)),\; as \; k\rightarrow\infty.
 $$


Taking $\psi\in C_0^\infty(0, T; C^\infty(D))$ as the testing
function for (\ref{ve1}) yields
\begin{eqnarray}\label{weakve}
&&-\int_0^T\hspace{-0.2cm}\big\langle \bar{v}_{\e_{j(k)}}(s),
\dot{\psi}(s) \big\rangle_Dds\\
&=&\hspace{-0.3cm}\int_0^T\hspace{-0.2cm}\big\langle
\bar{v}_{\e_{j(k)}}(s), \D\psi(s) \big\rangle_Dds
-\int_0^T\hspace{-0.2cm}\big\langle\p_\nu\psi(s),
\bar{v}_{\e_{j(k)}}(s) \big\rangle_{\Gamma_1}ds+\nonumber\\
&&\frac{1}{\sqrt{\e}}\int_0^T\hspace{-0.2cm}\big\langle
f(\bar{u}_{\e_{j(k)}}(s)-f(\bar{u}(s)), \psi\big\rangle_Dds
\hspace{-0.0cm}-\int_0^T\hspace{-0.2cm}\big\langle
\bar{v}_{\e_{j(k)}}(s), \psi(s)\big\rangle_{\Gamma_1}ds+ \nonumber
\\ &&\sqrt{\e}\int_0^T\hspace{-0.2cm}\big\langle \bar{u}_{\e_{j(k)}}(s),
\dot{\psi}(s)\big\rangle_{\Gamma_1}
ds+\int_0^T\hspace{-0.2cm}\big\langle
\sigma_2\dot{\overline{W}}_{2k}(s), \psi(s) \big\rangle_{\Gamma_1}ds
\nonumber.
\end{eqnarray}
We pass the limit $\e\rightarrow 0$ in (\ref{weakve}). Notice that
\begin{eqnarray*}
&&\Big| \sqrt{{\e_{j(k)}}}\int_0^T\big\langle
\bar{u}_{\e_{j(k)}}(s), \dot{\psi}(s)\big\rangle_{\Gamma_1} ds\Big|
\\ &\leq&
\sqrt{{\e_{j(k)}}}\int_0^T|\bar{u}_{\e_{j(k)}}(s)|_{L^2(\Gamma_1)}|\dot{\psi}(s)|_{L^2(\Gamma_1)}ds\\
&\leq&
\sqrt{{\e_{j(k)}}}\int_0^T|\bar{u}_{\e_{j(k)}}(s)|^2_{\mathbf{H}^\lambda_{\Gamma_2}(D)}ds+
\sqrt{{\e_{j(k)}}}\int_0^T|\dot{\psi}(s)|^2_{\mathbf{H}^1_{\Gamma_2}(D)}ds\\
&\rightarrow& 0,\;\;\e\rightarrow 0,\;\;{\rm
for}\;\;\omega\in\overline{\Omega}_\delta.
\end{eqnarray*}
By the assumption ($\mathbf{F}$) and (\ref{as})
$$
\frac{1}{\sqrt{{\e_{j(k)}}}}\big(f(\bar{u}_{\e_{j(k)}}(s))-f(\bar{u}(s))\big)
=f'(\tilde{u}_{\e_{j(k)}}(s))\bar{v}_{\e_{j(k)}}
$$
with $\tilde{u}_{\e_{j(k)}}\rightarrow \bar{u}$ almost surely on
$[0, T]\times D$. Then by Lemma \ref{Lions},
$f'(\tilde{u}_{\e_{j(k)}})$ converges weakly to $f'(\bar{u})$, for
$\omega\in\overline{\Omega}_\delta$. And by the choice of
$\bar{v}_{\e_{j(k)}}$, which converges strongly to $\bar{v}$ in
$L^2(0, T; L^2(D))$, we have
$f'(\tilde{u}_{\e_{j(k)}}(s))\bar{v}_{\e_{j(k)}}$ converges weakly
to $f'(\bar{u})\bar{v}$ in $L^2(0, T; L^2(D))$ which yields
$$
\frac{1}{\sqrt{\e_{j(k)}}}\int_0^T\big\langle
f(\bar{u}_{\e_{j(k)}}(s))-f(\bar{u}(s)),
\psi\big\rangle_Dds\rightarrow \int_0^T\langle f'(\bar{u})\bar{v},
\psi\rangle_Dds\;\;{\rm for}\;\;\omega\in\overline{\Omega}_\delta.
$$
Also by Lemma 3.1 of \cite{GyKry96}
$$
\int_0^T\big\langle \sigma_2\dot{\overline{W}}_{2k}(s), \psi(s)
\big\rangle_{\Gamma_1}ds\rightarrow \int_0^T\big\langle
\sigma_2\dot{\overline{W}}_{2}(s), \psi(s)
\big\rangle_{\Gamma_1}ds,\;\; k\rightarrow\infty
$$
in $\overline{\mathbb{P}}_\delta$-probability.

 Then combining all the above analysis for the terms in
(\ref{weakve}), we can pass the limit $\e\rightarrow 0$ in
(\ref{weakve}) and  conclude that
\begin{eqnarray}\label{weakv}
 &&\int_0^T\hspace{-0.2cm}\big\langle \dot{\bar{v}}(s), \psi(s) \big\rangle_Dds\\
 &=&\hspace{-0.3cm}\int_0^T\hspace{-0.2cm}\big\langle \D \bar{v}(s), \psi(s) \big\rangle_Dds
 +\int_0^T\hspace{-0.2cm}\big\langle f'(\bar{u})\bar{v}, \psi(s)\big\rangle_Dds
\nonumber\\ && \hspace{-0.2cm}-\int_0^T\hspace{-0.2cm}\big\langle
\bar{v}(s), \psi(s)
 \big\rangle_{\Gamma_1}ds+\int_0^T\hspace{-0.2cm}\big\langle \sigma_2\dot{\overline{W}}_2(s), \psi(s)
\big\rangle_{\Gamma_1}ds, \nonumber
\end{eqnarray}
which is the variational form of (\ref{v1})--(\ref{v4}). Notice that
we have proved the wellposedness of (\ref{v1})--(\ref{v4}).  Then by
the arbitrariness of $\delta$ and the same discussion in the proof
of Theorem \ref{main1}, we have the following result on normal
deviations.

\begin{theorem}\label{main2} (\textbf{Normal deviations principle}) \\
Assume that the conditions $(\mathbf{F})$ and $(\mathbf{\Sigma}')$
are satisfied. Let $u^0\in L^2(D)$ be a $(\mathcal{F}_0,
\mathcal{B}(L^2(D)))$- measurable random variable, which is
independent of $W(t)$, with
$\mathbf{E}|u^0|^{2}_{\mathbf{H}^1_{\Gamma_2}(D)}<\infty$ and
$\mathbf{E}\big[|u^0|^4_{L^2(D)}+ |\gamma_1 u^0|^4_{L^2(\Gamma_1)}
\big]<\infty$. Let $u_\e$ and $u$ be the unique weak solutions of
(\ref{fs1})--(\ref{fs5}) and (\ref{avsys1})--(\ref{avsys4}),
respectively. Then $\frac{1}{\sqrt{\e}}(u_\e-u)$ converges in
probability to a stochastic process v, which is the solution of the
linear random system (\ref{v1})--(\ref{v4}),  in the space $L^2(0,T;
L^2(D))$.
\end{theorem}

\vskip 0.8cm


\section{Large deviations  } \label{large}

In \S \ref{s5}, Theorem \ref{main1}, we have proved that
$u_\e\rightarrow u$
as $\e\rightarrow 0$. We have also obtained convergence result of
the normal deviations of order $\e^{\frac{1}{2}}$   in \S
\ref{s6}, Theorem \ref{main2}, which implies the normal deviations
of order $\e^\kappa$ tend to 0 for $0<\kappa<\frac{1}{2}$. In this
section we consider the logarithmic asymptotics of the deviations
of order $\e^\kappa$, $0<\kappa<\frac{1}{2}$, in probability. That
is, we consider the deviations of
$v^\kappa_\e=\e^{-\kappa}\left(u_\e-u \right)$ which satisfies
\begin{eqnarray}
\dot{v}^\kappa_\e&=&\D v^\kappa_\e+\e^{-\kappa}\big(f(u_\e)-f(u)\big),\;\; in\;\; D \label{vke1}\\
\p_\nu
v^\kappa_\e+v^\kappa_\e&=&-\e^{1-\kappa}\dot{u}_\e+\e^{\frac{1}{2}-\kappa}\sigma_2
\dot{W}_2,\;\; on
\;\;\Gamma_1 \label{vke2}\\
v^\kappa_\e&=&0,\;\; on \;\;\Gamma_2 \label{vke3} \\
v^\kappa_\e(0)&=&0. \label{vke4}
\end{eqnarray}
We intend to prove  that the family $\{v_\e^\kappa: \e>0\}$
satisfies the large deviations principle in $L^2(0, T; L^2(D))$.
We follow the results on large deviations in \cite{BD00} for
Polish space valued random elements; see also \cite{SS06} for
large deviations of two-dimensional
stochastic Navier-Stokes equations.\\

Let $H_0\subset H$ be Hilbert spaces with norm $|\cdot|_{H_0}$ and
$|\cdot|_H$ respectively. Assume that the embedding of $H_0$ in
$H$ is Hilbert-Schmidt. Define $\mathcal{A}$ the class of
$H_0$-valued $\mathcal{F}_t$-predictable process $w$ satisfying
$\int_0^T|w(s)|^2_{H_0}<\infty$ a.s. For $M>0$ let
$$
\mathcal{S}_M=\{w\in L^2(0, T; H_0):\int_0^T|w(s)|_{H_0}^2ds\leq M
\}
$$
which is a Polish space (i.e., complete separable metric space)
endowed with the weak topology. Define $\mathcal{A}_M=\{w\in
\mathcal{A}:w\in \mathcal{S}_M, a.s.\}$.

Let $E$ be a Polish space and $g^\e: C(0, T; H)\rightarrow E $ be
a measurable map. Let
  $V$ be an $H$-valued Wiener process. Define
$Y^\e=g^\e(V(\cdot))$. We consider the large deviation principle
for $Y^\e$ as $\e\rightarrow 0$. Since $E$ is a Polish space, the
Laplace principle and the large deviation principle are equivalent
\cite{SS06}.

\begin{defn}
A function $I$ mapping $E$ to $[0, \infty]$ is called  a rate
function if it is lower semicontinuous. A rate function $I$ is
called a good rate function if for each $M<\infty$, the level set
$\{y\in E: I(y)\leq M\}$ is compact in $E$.
\end{defn}
Recall that a family $\{Y^\e: \e>0\}$ of $E$-valued random
elements is said to satisfy the large deviations principle (LDP)
with speed $\alpha(\e)\rightarrow\infty$, as $\e\rightarrow 0$ and
rate function $I$ if (see \cite{FW84, PZ92})
\begin{enumerate}
  \item For any $\delta$, $\gamma>0$ and $y\in E$, there exists
  $\e_0>0$ such that for any $\e\in (0, \e_0)$
  $$
  \mathbf{P}\left\{|Y^\e-y|_E<\delta \right\}\geq \exp\{-\alpha(\e)(I(y)+\gamma)
  \}.
  $$
  \item For any $r$, $\delta$, $\gamma>0$, there exists $\e_0>0$
  such that for any $\e\in(0, \e_0)$
 $$
\mathbf{P}\left\{|Y^\e-\Phi(r)|_E\geq\delta \right\}\leq
\exp\{-\alpha(\e)(r-\gamma)\}
 $$
where $\Phi(r)=\{y\in E: I(y)\leq r\}$.
\end{enumerate}
It is well known that the large deviations principle and the
following Laplace principle are equivalent in Polish space.
\begin{defn}
Let $I$ be a rate function on metric space $E$. A family $\{Y^\e:
\e>0\}$ of $E$-valued random elements is said to satisfy the
Laplace principle on $E$ with rate function  $I$ and speed
$\alpha(\e)\rightarrow\infty$, $\e\rightarrow 0$, if for each real
valued, bounded and continuous function $\hbar$ defined on $E$,
$$
\lim_{\e\rightarrow 0}\frac{1}{\alpha(\e)}
\log\mathbf{E}\left\{\exp\left[-\alpha(\e)\hbar(Y^\e)\right]\right\}=-\inf_{y\in
E}\left\{\hbar(y)+I(y)\right\}.
$$
\end{defn}

For our purpose we introduce the assumption
\begin{description}
  \item[$(\mathbf{H})$] There exists a measurable map $g^0: C(0, T; H)\rightarrow
  E$ such that
   \begin{enumerate}
     \item Let $\{w^\e: \e>0\}\subset \mathcal{A}_M$ for some $M>0$.
     Let $w^\e$ converges in distribution to $w$. Then
     $g^\e\left(V(\cdot)+\frac{1}{\sqrt{\alpha(\e)}}\int_0^.w^\e(s)ds\right)$
     converges in distribution  to $g^0(\int_0^.w(s)ds)$.
     \item  For every $M<\infty$, the set $K_M=\{g^0(\int_0^.w(s)ds): w\in \mathcal{S}_M
     \}$ is a compact subset of $E$.
   \end{enumerate}

\end{description}
For each $g\in E$, define
\begin{equation}\label{rate}
I(y)=\inf_{\{w\in L^2(0, T;H_0):\;\; y=g^0(\int_0^.w(s)ds)
\}}\left\{\frac{1}{2}\int_0^T|w(s)|_{H_0}^2ds \right\}.
\end{equation}
Then we have the following theorem
\begin{theorem}\label{Lap}
Let $Y^\e=g^\e(V(\cdot))$. If $g^\e$ satisfies the assumption
$(\mathbf{H})$, then the family $\{Y^\e: \e>0\}$ satisfies the
Laplace principle in $E$ with rate function $I$ given by
(\ref{rate}) and speed $\alpha(\e)$.
\end{theorem}

The proof of the above theorem is   similar to   that of the proof
of Theorem 4.4 in \cite{BD00} which is for the speed
$\alpha(\e)=\e^{-1}$. We omit it here.

\vskip 0.5cm

In the following we apply the above result to the system
(\ref{vke1})--(\ref{vke4}). In this case $H=L^2(\Gamma_1)$,
$H_0=Q_2^{\frac{1}{2}}H$, $E=L^2(0, T; L^2(D))$,
$V(\cdot)=W_2(\cdot)$ and $\alpha(\e)=\e^{2\kappa-1}$,
$0<\kappa<\frac{1}{2}$. Since $Q_2$ is a trace class operator, the
embedding of $H_0$ in H is Hilbert-Schmidt. By the analysis of
Section \ref{s6} there exists a Borel measurable function $g^\e:
C(0, T; H)\rightarrow E$, such that $v^\kappa_\e=g^\e(W_2)$. We
intend to verify the assumption ($\mathbf{H}$) for $g^\e$. In fact
four lemmas are proved to complete the verification. Let $g^\e$ be
defined as above. For any $w\in\mathcal{A}_M$\,, $0<M<\infty$\,,
denote $g^\e\left(W_2(\cdot)+
\e^{-\frac{1}{2}+\kappa}\int_0^.w(s)ds\right)$ by
$\hat{v}_\e^\kappa$\,.

\begin{lemma}\label{H1}
 $\hat{v}_\e^\kappa$ is the unique weak solution of the following
stochastic system:
\begin{eqnarray}
\dot{\hat{v}}^\kappa_\e&=&\D \hat{v}^\kappa_\e+\e^{-\kappa}\big(f(u_\e)-f(u)\big),\;\; in\;\; D \label{hvke1}\\
\p_\nu
\hat{v}^\kappa_\e+\hat{v}^\kappa_\e&=&-\e^{1-\kappa}\dot{u}_\e+\sigma_2
w+\e^{\frac{1}{2}-\kappa}\sigma_2 \dot{W}_2,\;\; on
\;\;\Gamma_1 \label{hvke2}\\
\hat{v}^\kappa_\e&=&0,\;\; on \;\;\Gamma_2 \label{hvke3} \\
\hat{v}^\kappa_\e(0)&=&0 \label{hvke4}
\end{eqnarray}
in $L^2(\mathbf{\Omega}, C(0, T; L^2(D))\cap L^2(0, T;
\mathbf{H}^\lambda_{\Gamma_2}(D)))$. Here $\lambda$ is chosen in
Lemma \ref{H-lambda}.
\end{lemma}
\begin{proof}
This result   follows from a Girsanov argument. In fact let
$\widetilde{W}_2(\cdot)=W_2(\cdot)+\e^{-\frac{1}{2}+\kappa}\int_0^.w(s)ds$.
Then $\widetilde{W}_2$ is Wiener process with covariation $Q$
under the probability $\tilde{P}_w$ which satisfies
$$
d\tilde{P}_w=\exp\left\{
\e^{-\frac{1}{2}+\kappa}\int_0^Tw(s)dW_2(s)-\frac{1}{2}\e^{-1+2\kappa}\int_0^T|w(s)|_H^2ds
\right\}d\mathbf{P}.
$$
Then a similar analysis in Section \ref{s6} yields the result.
\end{proof}

\begin{lemma}\label{H2}
Let $w\in L^2(0, T; H)$. Then the following stochastic system
\begin{eqnarray}
\dot{\rho}_w&=&\D \rho_w+f'(u)\rho_w,\;\; in\;\; D \label{rho1}\\
\p_\nu \rho_w+\rho_w&=&\sigma_2w,\;\; on
\;\;\Gamma_1 \label{rho2}\\
\rho_w&=&0,\;\; on \;\;\Gamma_2 \label{rho3} \\
\rho_w(0)&=&0. \label{rho4}
\end{eqnarray}
has a unique weak solution $\rho_w\in C(0, T; L^2(D))\cap L^2(0, T;
\mathbf{H}^1_{\Gamma_2}(D))$.
\end{lemma}
\begin{proof}
This is a classical result of nonhomogeneous boundary problem
\cite{Lion1}.
\end{proof}

We now define the function $g^0$ as follows: $ g^0(h):=\rho_w$ if
$h=\int_0^.w(s)ds$ for some $w\in L^2(0,T; H)$, otherwise
$g^0(h)=0$. By the same discussion in Section \ref{s6} for the
normal deviations we   conclude that

\begin{lemma}\label{H3}
Let $\{w^\e\}\subset \mathcal{S}_M$ converge in distribution to
$w$, as a $\mathcal{S}_M$-valued random variable. Then
$g^\e\left(W_2(\cdot)+\e^{\kappa-\frac{1}{2}}\int_0^.w^\e(s)ds\right)$
converges in distribution to $g^0\left(\int_0^.w(s)ds\right )$ in
$E$.
\end{lemma}

\begin{lemma}\label{H4}
Let $0<M<\infty$ be fixed. Then the set $K_M=\{g^0(\int_0^.w(s)ds) :
w\in \mathcal{S}_M\}$ is compact in $E$.
\end{lemma}
\begin{proof}
By the definition of $\mathcal{S}_M$, for any sequence
$\{w^n\}\subset\mathcal{S}_M$, there is a subsequence $w^n$
(relabelled by n ) and $w\in \mathcal{S}_M$ such that $w^n$ weakly
converges to $w$ as $n\rightarrow\infty$. Then it is enough to
prove that $\rho_{w^n}$ converges to $\rho_w$ in $E$. Let
$\Theta_n=\rho_{w^n}-\rho_w$, then
\begin{eqnarray*}
\dot{\Theta}_n&=&\D \Theta_n+f'(u)\Theta_n\,,\;\; in\;\; D \\
\p_\nu \Theta_n+\Theta_n&=&\sigma_2(w^n-w)\,,\;\; on
\;\;\Gamma_1 \\
\Theta_n&=&0\,,\;\; on \;\;\Gamma_2 \\
\Theta_n(0)&=&0\,.
\end{eqnarray*}
A simple energy estimate  and the fact that embedding of $H_0$ in
$H$ is Hilbert-Schmidt yield
\begin{eqnarray*}
&&\sup_{0\leq t\leq
T}|\Theta_n(t)|^2_{L^2(D)}+\int_0^T|\Theta_n(s)|^2_{\mathbf{H}^1_{\Gamma_2}(D)}ds\\
&\leq& C(T, b, \sigma_2)\int_0^T|w^n(s)-w(s)|^2_{L^2(\Gamma_1)}ds\\
&\rightarrow & 0,\;\; as\;\; n\rightarrow\infty.
\end{eqnarray*}
This completes the proof.
\end{proof}

By the Lemma \ref{H1}--\ref{H4} and Theorem \ref{Lap} we can draw
the following result.

\begin{theorem}\label{main3} (\textbf{Large deviations principle}) \\
Assume that the conditions $(\mathbf{F})$ and $(\mathbf{\Sigma}')$
are satisfied. Let $u^0\in L^2(D)$ be a $(\mathcal{F}_0,
\mathcal{B}(L^2(D)))$-measurable random variable, which is
independent of $(W_1(t), W_2(t))$ with
$\mathbf{E}|u^0|^2_{\mathbf{H}^1_{\Gamma_2}(D)}<\infty$ and
$\mathbf{E}\big[|u^0|^4_{L^2(D)}+ |\gamma_1 u^0|^4_{L^2(\Gamma_1)}
\big]<\infty$. Let $u_\e$ and $u$ be the unique weak solutions of
(\ref{fs1})--(\ref{fs5}) and (\ref{sys1})--(\ref{sys4}),
respectively. Then for any $0<\kappa<\frac{1}{2}$,
$\e^{-\kappa}(u_\e-u)$ satisfies large deviations principle with
good rate function $I(\cdot)$ given by (\ref{rate}) and speed
$\e^{2\kappa-1}$ in $L^2(0, T; L^2(D))$.
\end{theorem}

\bigskip

{\bf Acknowledgements.}
   We would like to thank Francesco Russo, Mauro Mariani, Mark I. Freidlin,
    and Zhihui Yang for helpful discussions.
   We are grateful to Paul Dupuis for pointing out the paper
    \cite{BD00} and Padma Sundar for the paper \cite{SS06}.

\end{document}